\newcommand{\R}{\mathds R}
\newcommand{\boldR}{\mathds R}

\newcommand{\dd}{\mathrm d}
\newcommand{\Ddt}[1]{\tfrac{\mathrm D^{\hbox to 0pt{$\scriptscriptstyle#1$\hskip 0pt minus 1fil}}}{\dd t}}
\newcommand{\Id}{\mathrm{Id}}

\newcommand{\Lin}{\mathrm{Lin}}

\newcommand{\so}{\mathfrak{so}}
\newcommand{\ad}{\mathrm{ad}}

\newcommand{\tr}{\mathrm{tr}}
\newcommand{\gl}{\mathrm{gl}}

\newcommand{\Dim}{\mathrm{dim}}

\documentclass[oneside,draft,11pt]{amsart}

\usepackage{times}
\usepackage{dsfont}

\title[Characterization of Levi--Civita Connections]{Connections compatible with tensors.\\
A characterization of left-invariant Levi--Civita connections in Lie groups}
\author[P.\ Piccione]{Paolo Piccione}
\author[D.\ Tausk]{Daniel V.\ Tausk}
\address{Departamento de Matem\'atica,\hfill\break\indent  Universidade de S\~ao Paulo,
Brazil}
\email{piccione@ime.usp.br, tausk@ime.usp.br}
\urladdr{http://www.ime.usp.br/\~{}piccione,\hfill\break\phantom{URL: } http://www.ime.usp.br/\~{}tausk}

\date{September 27th, 2005}


\begin{document}


\theoremstyle{plain}\newtheorem{teo}{Theorem}[section]
\theoremstyle{plain}\newtheorem{prop}[teo]{Proposition}
\theoremstyle{plain}\newtheorem{lem}[teo]{Lemma}
\theoremstyle{plain}\newtheorem{cor}[teo]{Corollary}
\theoremstyle{definition}\newtheorem{defin}[teo]{Definition}
\theoremstyle{remark}\newtheorem{rem}[teo]{Remark}
\theoremstyle{plain} \newtheorem{assum}[teo]{Assumption}
\swapnumbers
\theoremstyle{definition}\newtheorem{example}{Example}[section]


\begin{abstract}
Symmetric connections that are compatible with semi-Riemannian metrics
can be characterized using an existence result for an integral leaf
of a (possibly non integrable) distribution. In this paper we give
necessary and sufficient conditions for a left-invariant connection
on a Lie group to be the Levi--Civita connection of some semi-Riemannian
metric on the group. As a special case, we will consider constant
connections in $\R^n$.
\end{abstract}

\maketitle

\begin{section}{Introduction}
In this short note we address the following problem: given a (symmetric) connection
$\nabla$ on a smooth manifold $M$, under which conditions there exists a semi-Riemannian
metric $g$ in $M$ which is $\nabla$-parallel?
This problem can be studied using holonomy theory (see \cite{KN}).
Alternatively, the problem can be cast in the language of distributions and integral submanifolds, as
follows. A connection $\nabla$ on a manifold $M$ induces naturally a connection in all
tensor bundles over $M$, in particular, on the bundle of all (symmetric) (2,0)-tensors
on $M$, say, $\nabla^{(2,0)}$. If $g$ is a $(2,0)$-tensor on $M$, $p\in M$ and $v,w\in T_pM$,
then the curvature $R^{(2,0)}(v,w)g$ is the bilinear form on $T_pM$ given by:
\begin{equation}\label{eq:Rbil}
\big(R^{(2,0)}(v,w)g\big)(\xi,\eta)=-g\big(R(v,w)\xi,\eta\big)-g\big(\xi,R(v,w)\eta\big),
\end{equation}
where  $\xi,\eta\in T_pM$  and $R$ is the curvature tensor of $\nabla$.
A semi-Riemannian metric is a (globally defined) symmetric nondegenerate $(2,0)$-tensor
on $M$, and
compatibility with $\nabla$ is equivalent to the property that the section is
everywhere tangent to the horizontal distribution determined by the
connection $\nabla^{(2,0)}$. However, such distribution is in general non integrable, namely,
integrability of the horizontal distribution is equivalent to the vanishing
of the curvature tensor $R^{(2,0)}$, which is equivalent to the vanishing of $R$.
Hence, the classical Frobenius theorem cannot be employed
in this situation. Nevertheless, the existence of simply \emph{one} integral
submanifold of a distribution, or, equivalently, of a parallel section
of a vector bundle endowed with a connection, may occur even in the case of non integrable distributions.
From \eqref{eq:Rbil}, one sees immediately that if $g$ is a
$(2,0)$-tensor on $M$, then the condition that $R^{(2,0)}$ vanishes
along $g$ is equivalent
to the condition of anti-symmetry of $g\big(R(v,w)\cdot,\cdot\big)$, for
all $p\in M$ and all $v,w\in T_pM$.

Let us consider the case that an open neighborhood  $V$ of a point $m_0$ of a manifold $M$
is \emph{ruled} by a family of curves issuing from $m_0$, parameterized by points of some manifold $\Lambda$.
What this means is that it is given a smooth function $\psi:Z\subset\R\times\Lambda\to M$,
defined on an open subset $Z$ of $\R\times\Lambda$, with $\psi(0,\lambda)=m_0$ for
all $\lambda$, and that admits a smooth
right inverse $\alpha:V\subset M\to Z$. Assume that it is given  a nondegenerate symmetric
bilinear form $g_0:T_{m_0}M\times T_{m_0}M\to\R$;
one obtains a semi-Riemannian metric on $V$ by \emph{spreading} $g_0$ with
parallel transport along the curves $t\mapsto\psi(t,\lambda)$.
If the tensor $g$ obtained in this way is such that $g\big(R(v,w)\cdot,\cdot\big)$
is an antisymmetric bilinear form on $T_pM$ for all $p\in V$ and all $v,w\in T_pM$,
then $g$ is $\nabla$-parallel. The precise statement of this fact is the following:
\begin{prop}\label{thm:vemdemetrica}
Let $M$ be a smooth manifold, $\nabla$ be a symmetric connection on $TM$, $m_0\in M$ and $g_0$ be a nondegenerate
symmetric bilinear form on $T_{m_0}M$. Let $\psi:Z\subset\R\times\Lambda\to M$ be
a $\Lambda$-parametric family of curves on $M$ with a local right inverse $\alpha:V\subset M\to Z$; assume
that $\psi(0,\lambda)=m_0$, for all $\lambda\in M$. For each $(t,\lambda)\in Z$, we denote by
$P_{(t,\lambda)}:T_{m_0}M\to T_{\psi(t,\lambda)}M$ the parallel transport along $t\mapsto\psi(t,\lambda)$.
Assume that for all $(t,\lambda)\in Z$ the linear operator:
\begin{equation}\label{eq:hipoteseR}
P_{(t,\lambda)}^{-1}\big[R_{\psi(t,\lambda)}(v,w)\big]P_{(t,\lambda)}:T_{m_0}M\longrightarrow T_{m_0}M
\end{equation}
is anti-symmetric with respect to $g_0$, for all $v,w\in T_{\psi(t,\lambda)}M$,
where
\[R_{\psi(t,\lambda)}(v,w):T_{\psi(t,\lambda)}M\longrightarrow T_{\psi(t,\lambda)}M\]
denotes the linear operator corresponding to the curvature tensor of $\nabla$. Then $\nabla$ is the Levi-Civita
connection of the semi-Riemannian metric $g$ on $V\subset M$ defined by setting:
\[g_m(\cdot,\cdot)=g_0(P_{\alpha(m)}^{-1}\cdot,P_{\alpha(m)}^{-1}\cdot),\]
for all $m\in V$.
\end{prop}
\begin{proof}
See \cite{Frobenius}
\end{proof}
In the real analytic case, we have the following global result:
\begin{prop}\label{thm:propmetricarealanalitica}
Let $M$ be a  simply-connected real-analytic manifold and let $\nabla$ be a real-analytic
symmetric connection on $TM$. If there exists a semi-Riemannian metric $g$ on a nonempty open connected subset
of $M$ having $\nabla$ as its Levi-Civita connection then $g$ extends to a globally defined semi-Riemannian metric
on $M$ having $\nabla$ as its Levi-Civita connection.\qed
\end{prop}
The two results above will be used in Sections~\ref{sec:constconnRn}, \ref{sec:twodimcase}
and \ref{sec:leftinvconn}
to characterize symmetric connections in Lie groups that are constant in left invariant referentials.
The case of $\R^n$ (Lemma~\ref{thm:quasela} and Proposition~\ref{thm:sog0}), and more specifically the $2$-dimensional case
(Proposition~\ref{thm:propcasodim2}), will be studied with
some more detail.
\smallskip

It is an interesting problem to study conditions for the existence, uniqueness,
multiplicity, etc., of (symmetric) connections that are compatible with arbitrarily
given tensors. It is well known that semi-Riemannian metrics admit \emph{exactly one}
symmetric and compatible connection, called the \emph{Levi--Civita connection} of the metric.
Uniqueness can be deduced also by a curious combinatorial argument, see Corollary~\ref{thm:curious}.
The next interesting case is that of symplectic forms, in which case one has existence, but
not uniqueness. We will start the paper with a short section containing  a couple
of simple results concerning compatible connections. First,  we will show the combinatorial
argument that shows the uniqueness of the Levi--Civita connection of a semi-Riemannian
metric tensor (Corollary~\ref{thm:curious}). Second, we will prove that the
existence of a symmetric connection compatible with a nondegenerate two-form $\omega$ is
equivalent to the fact that $\omega$ is closed, in which case there are infinitely many
symmetric connections compatible with $\omega$ (Lemma~\ref{thm:exmultsymconsympl}).
\end{section}

\begin{section}{Connections compatible with tensors}
Let $M$ be a smooth manifold and let $\tau$ be any tensor in $M$; we will be mostly interested in the case
when $\tau=g$ is a \emph{semi-Riemannian metric tensor} on $M$ (i.e., $\tau$ is a nondegenerate symmetric $(2,0)$-tensor),
or when $\tau=\omega$ is a \emph{symplectic form} on $M$ (i.e., $\tau$ is a nondegenerate closed $2$-form).
If $\nabla$ is a connection in $M$, i.e., a connection on the tangent bundle $TM$, then we have naturally
induced connections on all tensor bundles on $M$, all of which will be denoted by the same symbol $\nabla$.

The \emph{torsion} of $\nabla$ is the anti-symmetric tensor
\[T(X,Y)=\nabla_XY-\nabla_YX-[X,Y],\]
where $[X,Y]$ denotes the Lie brackets of the vector fields $X$ and $Y$; $\nabla$ is called \emph{symmetric}
if $T=0$.
The connection $\nabla$ is said to be \emph{compatible} with $\tau$ if $\tau$ is $\nabla$-parallel, i.e.,
when $\nabla\tau=0$.

Establishing whether a given tensor $\tau$ admits compatible
connections is a \emph{local} problem. Namely, one can use
partition of unity to extend locally defined connections and
observe that a convex combination of compatible connections is a
compatible connection. In local coordinates, finding a connection
compatible with a given tensor reduces to determining the
existence of solutions for a non homogeneous linear system
for the Christoffel symbols of the connection.

It is well known that semi-Riemannian
metric tensors admit a \emph{unique} compatible symmetric
connection, called the Levi--Civita connection of the metric
tensor, which can be given explicitly by Koszul formula (see for instance \cite{dC}).
Uniqueness of the Levi--Civita connection can be obtained by a curious
combinatorial argument, as follows.

Suppose  that $\nabla$ and
$\widetilde\nabla$ are connections on $M$; their difference $\widetilde\nabla-\nabla$ is
a \emph{tensor}, that will be denoted by $\mathfrak t$:
\[\mathfrak t(X,Y)=\widetilde\nabla_XY-\nabla_XY,\]
where $X$ and $Y$ are smooth vector fields on $M$.
If both $\nabla$ and $\widetilde\nabla$ are symmetric connections, then
$\mathfrak t$ is symmetric:
\[\mathfrak t(X,Y)-\mathfrak t(Y,X)=\widetilde\nabla_XY-\nabla_XY-\widetilde\nabla_YX+\nabla_YX=[X,Y]+[Y,X]=0.\]

\begin{lem}\label{thm:symantisym} Let $U$ be a set and $\rho:U\times U\times
U\to\R$ be a map that is symmetric in its first two variables and
anti-symmetric in its last two variables. Then $\rho$ is
identically zero.
\end{lem}
\begin{proof}
Let $u_1,u_2,u_3\in U$ be fixed. We have:
$$\rho(u_1,u_2,u_3)=\rho(u_2,u_1,u_3)=-\rho(u_2,u_3,u_1)=-\rho(u_3,u_2,u_1),$$
so that $\rho$ is anti-symmetric in the first and the third
variables. On the other hand:
$$\rho(u_1,u_2,u_3)=-\rho(u_3,u_2,u_1)=-\rho(u_2,u_3,u_1)=\rho(u_1,u_3,u_2),$$
so that $\rho$ is symmetric in the second and the third variables.
This concludes the proof.
\end{proof}

\begin{cor}\label{thm:curious}
There exists at most one symmetric connection which is compatible with a semi-Riemannian
metric.
\end{cor}
\begin{proof}
Assume that $g$ is a semi-Riemannian metric on $M$, and let $\nabla$ and $\widetilde\nabla$ two
symmetric connections such that $\nabla g=\widetilde\nabla g=0$; for all $p\in M$
consider the map $\rho:T_pM\times T_pM\times T_pM\to\R$ given by:
\[\rho(X,Y,Z)=g\big(\mathfrak t(X,Y),Z\big),\]
where $\mathfrak t$ is the difference $\widetilde\nabla-\nabla$.
Since $\mathfrak t$ is symmetric, then $\rho$ is symmetric in the first two variables.
On the other hand, $\rho$ is anti-symmetric in the last two variables:
\begin{multline*}\rho(X,Y,Z)+\rho(X,Z,Y)\\=g\big(\widetilde\nabla_XY,Z\big)-g\big(\nabla_XY,Z\big)+
g\big(\widetilde\nabla_XZ,Y\big)-g\big(\nabla_XZ,Y\big)\\=\widetilde\nabla g(X,Y,Z)-\nabla g(X,Y,Z)=0.
\end{multline*}
By Lemma~\ref{thm:symantisym}, $\rho=0$, hence $\mathfrak t=0$, and thus $\widetilde\nabla=\nabla$.
\end{proof}
For symplectic forms, the situation changes radically.
Among all nondegenerate two-forms, the existence of a symmetric compatible
connection characterizes the symplectic ones:
\begin{lem}\label{thm:exmultsymconsympl}
Let $\omega$ be a nondegenerate $2$-form on a (necessarily even dimensional) manifold $M$.
There exists a symmetric connection in $M$ compatible
with $\omega$ if and only if $\omega$ is closed. In this case, there are infinitely
many symmetric connections that are compatible with $\omega$.
\end{lem}
\begin{proof}
If $\omega$ is closed, i.e., if $\omega$ is a symplectic form on
$M$, Darboux theorem tells us that one can find coordinates
$(q,p)$ around every point of $M$ such that
$\omega=\sum_{i}\mathrm dq^i\wedge\mathrm dp_i$, which means that
$\omega$ is constant in such coordinate system. The (locally
defined) symmetric connection which has vanishing Christoffel
symbols in such coordinates is clearly compatible with $\omega$.
As observed above, using partitions of unity one can find a
globally defined symmetric connection compatible with $\omega$.

\noindent Conversely, if $\nabla$ is any symmetric   connection in
$M$, then $\mathrm d\omega$ is given by
$\frac12\mathrm{Alt}(\nabla\omega)$, where $\mathrm{Alt}$ denotes
the alternator; in particular, if there exists a compatible
symmetric connection it must be $\mathrm d\omega=0$.
\end{proof}

\end{section}

\begin{section}{Constant connections in $\mathbf{\boldR^n}$}
\label{sec:constconnRn}

Let $\Gamma:\R^n\times\R^n\to\R^n$ be a symmetric bilinear map and consider the symmetric connection
$\nabla$ on $\R^n$ defined by:
\begin{equation}\label{eq:nablaGamma}
\nabla_XY=\dd Y(X)+\Gamma(X,Y),
\end{equation}
for any smooth vector fields $X$, $Y$ on $\R^n$. We now apply the result of Proposition~\ref{thm:vemdemetrica}
to determine when $\nabla$ is the Levi-Civita connection of a semi-Riemannian metric on $\R^n$.
Given $v\in\R^n$ then the parallel transport $P_{(t,v)}:\R^n\to\R^n$ along the curve $t\mapsto tv$ is given by:
\[P_{(t,v)}=\exp\big({-t\,\Gamma(v)}\big),\]
where we identify $\Gamma$ with the linear map $\R^n\ni v\mapsto\Gamma(v,\cdot)\in\Lin(\R^n)$.
For any $v,w\in\R^n$, the curvature tensor $R$ of $\nabla$ is given by:
\[R_x(v,w)=\Gamma(v)\Gamma(w)-\Gamma(w)\Gamma(v)=[\Gamma(v),\Gamma(w)]\in\Lin(\R^n),\]
for all $x\in\R^n$. Applying Proposition~\ref{thm:vemdemetrica} to the $\R^n$-parametric family
of curves $\psi(t,\lambda)=t\lambda\in\R^n$ with right inverse $\alpha:\R^n\ni v\mapsto(1,v)\in\R\times\R^n$
we obtain the following:
\begin{lem}\label{thm:quasela}
Let $g_0$ be a nondegenerate symmetric bilinear form on $T_0\R^n\cong\R^n$. Then $g_0$ extends to a semi-Riemannian
metric $g$ on $\R^n$ having \eqref{eq:nablaGamma} as its Levi-Civita connection if and only if
the linear operator:
\[\exp\big(\Gamma(v)\big)[\Gamma(w_1),\Gamma(w_2)]\exp\big({-\Gamma(v)}\big):\R^n\to\R^n\]
is anti-symmetric with respect to $g_0$, for all $v,w_1,w_2\in\R^n$.\qed
\end{lem}

Given a nondegenerate symmetric bilinear form $g_0$ on $\R^n$ we denote by $\so(g_0)$ the Lie algebra
of all $g_0$-anti-symmetric endomorphisms of $\R^n$. Given a linear endomorphism $X$ of $\R^n$ we write:
\[\ad_X(Y)=[X,Y]=XY-YX,\]
for all $Y\in\Lin(\R^n)$.
\begin{prop}\label{thm:sog0}
Let $\Gamma:\R^n\times\R^n\to\R^n$ be a symmetric bilinear map and let $\mathcal S\subset\Lin(\R^n)$ be the range
of the linear map $v\mapsto\Gamma(v,\cdot)$. A nondegenerate symmetric bilinear form $g_0$ on $T_0\R^n\cong\R^n$
extends to a semi-Riemannian metric $g$ on $\R^n$ having \eqref{eq:nablaGamma} as its Levi-Civita connection
if and only if:
\[(\ad_X)^k[Y,Z]\in\so(g_0),\]
for all $X,Y,Z\in\mathcal S$ and all $k\ge0$.
\end{prop}
\begin{proof}
By Lemma~\ref{thm:quasela}, $g_0$ extends to a semi-Riemannian metric $g$ on $\R^n$ having \eqref{eq:nablaGamma} as
its Levi-Civita connection if and only if:
\[\exp(tX)[Y,Z]\exp({-tX})\in\so(g_0),\]
for all $X,Y,Z\in\mathcal S$ and all $t\in\R$. The conclusion follows by observing that:
\[\exp(tX)[Y,Z]\exp({-tX})=\sum_{k=0}^{+\infty}\frac{t^k}{k!}(\ad_X)^k[Y,Z].\qedhere\]
\end{proof}

\begin{cor}\label{thm:corgg}
Let $\Gamma:\R^n\times\R^n\to\R^n$ be a symmetric bilinear map and let $\mathcal S\subset\Lin(\R^n)$ be the range
of the linear map $v\mapsto\Gamma(v,\cdot)$. Denote by $\mathfrak g$ the Lie algebra spanned by $\mathcal S$ and by
$\mathfrak g'=[\mathfrak g,\mathfrak g]$ the commutator subalgebra of $\mathfrak g$. If $\mathfrak g'$
is contained in $\so(g_0)$ for some nondegenerate symmetric bilinear form $g_0$ on $T_0\R^n\cong\R^n$
then $g_0$ extends to a semi-Riemannian metric $g$ on $\R^n$ having \eqref{eq:nablaGamma} as its Levi-Civita connection.\qed
\end{cor}

\end{section}

\begin{section}{The two-dimensional case}
\label{sec:twodimcase}
If $n=2$, the Lie algebra $\so(g_0)$ is one-di\-mensional. This observation allows us to show that, for $n=2$,
the condition $\mathfrak g'\subset\so(g_0)$ in the statement
of Corollary~\ref{thm:corgg} is also necessary for $g_0$ to extend to a semi-Riemannian metric $g$
on $\R^2$ having \eqref{eq:nablaGamma} as its Levi-Civita connection.
\begin{lem}\label{thm:cond2x2}
Let $\Gamma:\R^2\times\R^2\to\R^2$ be a symmetric bilinear map and let $\mathcal S\subset\Lin(\R^2)$ be the range
of the linear map $v\mapsto\Gamma(v,\cdot)$. Denote by $\mathfrak g$ the Lie algebra spanned by $\mathcal S$.
Then a nondegenerate symmetric bilinear form $g_0$ on $T_0\R^2\cong\R^2$ extends to a semi-Riemannian metric $g$ on $\R^2$
having \eqref{eq:nablaGamma} as its Levi-Civita connection if and only if $\mathfrak g'\subset\so(g_0)$.
\end{lem}
\begin{proof}
Define a sequence $\mathcal S_k$ of subspaces of $\mathfrak g$ inductively by setting $\mathcal S_1=\mathcal S$
and by taking $\mathcal S_{k+1}=[\mathcal S,\mathcal S_k]$ to be the linear span of all commutators $[X,Y]$, with
$X\in\mathcal S$, $Y\in\mathcal S_k$. Using the Jacobi identity it is easy to show that $[\mathcal S_k,\mathcal S_l]\subset
\mathcal S_{k+l}$ and therefore:
\[\mathfrak g=\sum_{k=1}^\infty\mathcal S_k,\quad\mathfrak g'=\sum_{k=2}^\infty \mathcal S_k.\]
By Proposition~\ref{thm:sog0}, if extends to a semi-Riemannian metric $g$ on $\R^2$
having \eqref{eq:nablaGamma} as its Levi-Civita connection then $\mathcal S_2$ and $\mathcal S_3$ are contained
in $\so(g_0)$. Since $\so(g_0)$ is one dimensional, we have either $\mathcal S_3=0$ or $\mathcal S_3=\mathcal S_2$;
in the first case, $\mathcal S_k=0$ for all $k\ge3$ and in the latter case $\mathcal S_k=\mathcal S_2$ for all
$k\ge3$. In any case, $\mathfrak g'=\mathcal S_2$ and the conclusion follows.
\end{proof}

\begin{lem}\label{thm:cond2x22}
Let $A:\R^2\to\R^2$ be a nonzero linear map. There exists a nondegenerate symmetric bilinear form $g_0$
on $\R^2$ with $A\in\so(g_0)$ if and only if $\tr\,A=0$ and $\det A\ne0$; moreover, $g_0$ is positive definite
(resp., has index $1$) if and only if $\det A>0$ (resp., $\det A<0$).
\end{lem}
\begin{proof}
Assume that $\tr\,A=0$ and $\det A\ne0$. Write $\det A=-\epsilon a^2$, with $\epsilon=\pm1$ and $a>0$. It is easy to see that $A$ is represented by the
matrix $\left(\begin{smallmatrix}0&\epsilon a\\a&0\end{smallmatrix}\right)$ in some basis $(b_1,b_2)$ of $\R^2$.
We define $g_0$ by setting:
\begin{equation}\label{eq:defg0}
g_0(b_1,b_2)=0,\quad g_0(b_1,b_1)=1,\quad g_0(b_2,b_2)=-\epsilon.
\end{equation}
Conversely, if $A\in\so(g_0)$ for some $g_0$ then we can choose a basis $(b_1,b_2)$ of $\R^2$ such that
\eqref{eq:defg0} holds and the matrix of $A$ on such basis is of the form
$\left(\begin{smallmatrix}0&\epsilon a\\a&0\end{smallmatrix}\right)$.
\end{proof}

\begin{cor}\label{thm:corglinha}
Let $\Gamma:\R^2\times\R^2\to\R^2$ be a symmetric bilinear map and let $\mathcal S\subset\Lin(\R^2)$ be the range
of the linear map $v\mapsto\Gamma(v,\cdot)$. Denote by $\mathfrak g$ the Lie algebra spanned by $\mathcal S$.
There exists a semi-Riemannian metric on $\R^2$ having \eqref{eq:nablaGamma} as its Levi-Civita connection if and only if
either $\mathfrak g'=0$ or $\mathfrak g'$ is one-dimensional and it is spanned by an invertible $2\times2$ matrix.
\end{cor}
\begin{proof}
Follows from Lemmas~\ref{thm:cond2x2} and \ref{thm:cond2x22}, observing that the elements of $\mathfrak g'$ have
null trace.
\end{proof}

\begin{lem}\label{thm:temcentro}
Let $\mathfrak g$ be a three-dimensional real Lie algebra with $\mathfrak g'$ one-dimen\-sional. Then the center
$\mathfrak z(\mathfrak g)$ of $\mathfrak g$ is one-dimensional.
\end{lem}
\begin{proof}
Let $Z$ denote a generator of $\mathfrak g'$, so that $[X,Y]=\alpha(X,Y)Z$, for all $X,Y\in\mathfrak g$, where
$\alpha$ is an antisymmetric bilinear form on $\mathfrak g$; clearly, the kernel of $\alpha$ is the center of $\mathfrak g$.
Since $\mathfrak g$ is three-dimensional, the kernel
of $\alpha$ is either $\mathfrak g$ or it is one-dimensional; the first possibility does not occur, since $\mathfrak g'$
is nonzero.
\end{proof}

\begin{cor}\label{thm:possibilities}
Let $\mathfrak g$ be a three-dimensional real Lie algebra with $\mathfrak g'$ one-dimen\-sional. Then there
exists a basis $(X,Y,Z)$ of $\mathfrak g$ such that one the following commutation relations holds:
\begin{enumerate}
\item\label{itm:rel1} $[X,Y]=[X,Z]=0$, $[Y,Z]=X$;
\item\label{itm:rel2} $[X,Y]=[Y,Z]=0$, $[Z,X]=X$.
\end{enumerate}
\end{cor}
\begin{proof}
Choose a basis $(X,Y,Z)$ of $\mathfrak g$ with $X$ in $\mathfrak g'$. If $\mathfrak z(\mathfrak g)=\mathfrak g'$
then $[Y,Z]\ne0$, otherwise $\mathfrak g'=0$; thus, we can replace $Z$ with a scalar multiple of $Z$ so that
$[Y,Z]=X$ and relations \eqref{itm:rel1} hold. If $\mathfrak z(\mathfrak g)\ne\mathfrak g'$, we may assume
that $Y\in\mathfrak z(\mathfrak g)$ and $[Z,X]\ne0$; again, replacing $Z$ with a scalar multiple of $Z$
gives $[Z,X]=X$ and relations \eqref{itm:rel2} hold.
\end{proof}

In what follows we denote by $\gl(n,\R)$ the Lie algebra of linear endomorphisms of $\R^n$.
\begin{lem}\label{thm:temZeX}
Let $\mathfrak g$ be a three-dimensional Lie subalgebra of $\gl(2,\R)$ with $\mathfrak g'$ one-dimensional.
There exists a basis $(X,Y,Z)$ of $\mathfrak g$ with $Y=\Id$ and $[Z,X]=X$.
\end{lem}
\begin{proof}
We show that $\Id$ is in $\mathfrak g$. Assume not. By Lemma~\ref{thm:temcentro},
there exists a nonzero element $W$ in $\mathfrak z(\mathfrak g)$. Then $W$ commutes with $\mathfrak g$ and
with $\Id$, which implies that $W$ is in the center of $\gl(2,\R)$; thus $W$ is a nonzero multiple
of $\Id$, contradicting our assumption.

Now $\Id\in\mathfrak g$ implies that $\Id$ spans $\mathfrak z(\mathfrak g)$; thus possibility
\eqref{itm:rel1} in the statement of Corollary~\ref{thm:possibilities} does not occur for it would imply
that $[Y,Z]$ is a nonzero multiple of the identity. Hence possibility \eqref{itm:rel2} occurs and we can
assume that $Y=\Id$.
\end{proof}

\begin{lem}\label{thm:twodimeng}
If $\mathfrak g$ is a two-dimensional real Lie algebra with $\mathfrak g'\ne0$ then there exists a basis
$(X,Z)$ of $\mathfrak g$ with $[Z,X]=X$.
\end{lem}
\begin{proof}
Let $X$ be a nonzero element in $\mathfrak g'$; clearly, $\mathfrak g'$ is one-dimensional. We can choose
$Z\in\mathfrak g$, $Z\not\in\mathfrak g'$, with $[Z,X]=X$.
\end{proof}

\begin{lem}\label{thm:Xnotinvertible}
If $X,Z\in\gl(n,\R)$ and $[Z,X]=X$ then $X$ is not invertible.
\end{lem}
\begin{proof}
If $X$ were invertible then $[Z,X]=X$ would imply $Z-XZX^{-1}=\Id$. A contradiction
is obtained by taking traces on both sides.
\end{proof}

\begin{prop}\label{thm:propcasodim2}
Let $\Gamma:\R^2\times\R^2\to\R^2$ be a symmetric bilinear map and let $\mathcal S\subset\Lin(\R^2)$ be the range
of the linear map $v\mapsto\Gamma(v,\cdot)$. Then there exists a semi-Riemannian metric $g$ on $\R^2$ having
\eqref{eq:nablaGamma} as its Levi-Civita connection if and only if $[X,Y]=0$, for all $X,Y\in\mathcal S$.
In this case, a semi-Riemannian metric $g$ on $\R^2$ having \eqref{eq:nablaGamma} as its Levi-Civita connection
can be chosen with an arbitrary value $g_0$ at the origin.
\end{prop}
\begin{proof}
If $[X,Y]=0$ for all $X,Y\in\mathcal S$ then, by Lemma~\ref{thm:cond2x2}, any nondegenerate symmetric
bilinear form $g_0$ on $T_0\R^2\cong\R^2$ extends to semi-Riemannian metric $g$ on $\R^2$ having
\eqref{eq:nablaGamma} as its Levi-Civita connection. Now assume that there exists a semi-Riemannian metric
on $\R^2$ having \eqref{eq:nablaGamma} as its Levi-Civita connection and denote by $\mathfrak g$ the Lie
algebra spanned by $\mathcal S$. By Corollary~\ref{thm:corglinha}, either $\mathfrak g'=0$ or $\mathfrak g'$
is one-dimensional and it is spanned by an invertible $2\times2$ matrix. Let us show that the second possibility
cannot occur. If $\mathfrak g'\ne0$ then $2\le\Dim(\mathfrak g)\le4$. If $\Dim(\mathfrak g)=4$ then $\mathfrak g=\gl(2,\R)$
and $\mathfrak g'$ is three-dimensional, which is not possible. If either $\Dim(\mathfrak g)=2$ or $\Dim(\mathfrak g)=3$
then by Lemmas~\ref{thm:temZeX} and \ref{thm:twodimeng} there exist $X,Z\in\mathfrak g$ with $[Z,X]=X$ and
such that $X$ spans $\mathfrak g'$. By Lemma~\ref{thm:Xnotinvertible}, $X$ is not invertible and we obtain
a contradiction.
\end{proof}

\end{section}

\begin{section}{Left-Invariant Connections on Lie Groups}
\label{sec:leftinvconn}
Let $G$ be a Lie group and $\nabla$ be a left-invariant connection on $G$. The connection $\nabla$ is determined
by a bilinear map $\Gamma:\mathfrak g\times\mathfrak g\to\mathfrak g$, i.e.:
\[\nabla_XY=\Gamma(X,Y),\]
for any left-invariant vector fields $X$, $Y$ on $G$.

The torsion of $\nabla$ is given by:
\[T(X,Y)=\Gamma(X,Y)-\Gamma(Y,X)-[X,Y].\]
Observe that $\nabla$ is torsion-free if and only if there exists a symmetric bilinear map $B:\mathfrak g\times\mathfrak g\to\mathfrak g$
with $\Gamma(X,Y)=B(X,Y)+\frac12[X,Y]$, for all $X,Y\in\mathfrak g$. If we identify $\Gamma$ with the linear map
$\mathfrak g\ni X\mapsto\Gamma(X,\cdot)\in\Lin(\mathfrak g)$ then $\nabla$ is torsion-free if and only if $\Gamma$ is
a Lie algebra homomorphism. The curvature tensor of $\nabla$ is given by:
\[R(X,Y)=[\Gamma(X),\Gamma(Y)]-\Gamma\big([X,Y]\big);\]
observe that the first bracket is the commutator in $\Lin(\mathfrak g)$ and the second is the Lie algebra product
of $\mathfrak g$.

Given a curve $\gamma$ on $G$, we identify vector fields along $\gamma$ with curves on $\mathfrak g$ by
left translation. Using this identification, the parallel transport of $Y\in\mathfrak g$ along a one-parameter subgroup
$t\mapsto\exp(tX)\in G$ is given by $t\mapsto e^{-t\Gamma(X)}Y\in\mathfrak g$.

\begin{prop}\label{thm:localnablagroup}
Assume that $\nabla$ is torsion-free and let $h:\mathfrak g\times\mathfrak g\to\R$ be a nondegenerate symmetric bilinear
form. The following condition
is necessary and sufficient for the existence of an extension of $h$ to a semi-Riemannian metric on a neighborhood
of the identity of $G$ whose Levi-Civita connection is $\nabla$:
\begin{equation}\label{eq:condGlocal}
e^{\Gamma(Z)}\Big([\Gamma(X),\Gamma(Y)]-\Gamma\big([X,Y]\big)\Big)e^{-\Gamma(Z)}\in\so(h),\quad\text{for all $X,Y,Z\in\mathfrak g$}.
\end{equation}
In \eqref{eq:condGlocal} we have denoted by $\so(h)$ the Lie subalgebra of $\Lin(\mathfrak g)$ consisting of $h$-anti-symmetric
linear operators.
\end{prop}
\begin{proof}
Set $\Lambda=\mathfrak g$ and consider the one-parameter family of curves $\psi:\R\times\Lambda\to G$ defined by
$\psi(t,\lambda)=\exp(t\lambda)$. If $U$ is an open neighborhood of the origin of $\mathfrak g$ that is mapped
diffeomorphically by $\exp$ onto an open neighborhood $\exp(U)$ of the identity of $G$ then a local right inverse
for $\psi$ can be defined by setting $\alpha(g)=\big(1,(\exp\vert_U)^{-1}(g)\big)$, for all $g\in\exp(U)$.
The conclusion follows from Proposition~\ref{thm:vemdemetrica}.
\end{proof}

\begin{cor}
Assume that $\nabla$ is torsion-free and let $h:\mathfrak g\times\mathfrak g\to\R$ be a nondegenerate symmetric bilinear
form. If $G$ is (connected and) simply-connected then condition \eqref{eq:condGlocal} is necessary and sufficient
for the existence of a globally defined semi-Riemannian metric on $G$ whose Levi-Civita connection is $\nabla$.
\end{cor}
\begin{proof}
It follows from Proposition~\ref{thm:localnablagroup} and from Proposition~\ref{thm:propmetricarealanalitica}
observing that left-invariant objects on a Lie group are always real-analytic.
\end{proof}

\begin{lem}
Condition \eqref{eq:condGlocal} is equivalent to:
\[\ad_{\Gamma(Z)}^n\Big([\Gamma(X),\Gamma(Y)]-\Gamma\big([X,Y]\big)\Big)\in\so(h),\quad\text{for all $X,Y,Z\in\mathfrak g$,
$n\ge0$},\]
where $\ad_A(B)=[A,B]$, for all $A,B\in\Lin(\mathfrak g)$.
\end{lem}
\begin{proof}
Replace $Z$ by $tZ$ in \eqref{eq:condGlocal} and compute the Taylor expansion in powers of $t$ of the corresponding
expression.
\end{proof}


\end{section}

\end{document}